\documentclass[reqno,dvips,12pt]{article}
\usepackage{epsfig,graphicx,amsmath,amssymb,amsfonts}
cm \textwidth=16.0 true cm \emergencystretch=10pt

\newtheorem{theorem}{\qquad Theorem}
\newtheorem{lemma}{\qquad Lemma} %[section]

\DeclareMathOperator{\re}{re}

\newcommand{\cA}{{\cal A}}
\newcommand{\cM}{{\cal M}}
\newcommand{\cE}{{\cal E}}

\newcommand{\cN}{{\cal N}}
\newcommand{\cP}{{\cal P}}

\newcommand{\nn}{\nonumber}

\newcommand{\pa}{\partial}

\newcommand{\e}{\varepsilon}

\newcommand{\abs}[1]{\left\vert{#1}\right\vert}

\begin{document}

\title {\textbf{Rank Distributions in Semiotics}}

\author {V.P.Maslov, T.V.Maslova\thanks{Moscow Institute of Electronics and
Mathematics, Moscow Institute of Economics; pm@miem.edu.ru:}}
\date{}
\maketitle

\begin{abstract}
The notions of real and user cardinality of a sign are introduced.
Rank distributions can be extended to arbitrary sign objects,
i.e., semiotic systems. The dynamics of the distribution of
consumer durables, such as automobiles, is studied.
\end{abstract}

Usually, in semiotics only relatively short strings of signs of
signs (discourses) are considered, while long strings with large
parameters have not really been studied. We shall introduce
notions generalizing the notion of participant in communicative
and non-communicative sign systems: instead of the terms
"narrateur" and "narrataire," or "interlocuteur" and
"interlocutaire" (see \cite{Semiotica}, p. 508]), we shall use
the pair of terms \textit{generator} and \textit{user}.

Semiotic objects, i.e., signs, can be of different types
\cite{Semiotica,Koma}.

A word of a natural language is \textit{a sign}. The collection of
words is \textit{the dictionary of signs}. We use the term
dictionary of signs rather than ''alphabet of signs'' to stress
that the number of signs can be very large. The activity index of
a sign is the number of its occurrences. We shall call this index
the real cardinality of the sign. \textit{The real cardinality of
the dictionary of signs} is the total number of occurrences of all
the signs from the dictionary (collection of signs). In language
systems, the cardinality of a dictionary (collection of words)
corresponds to the number of occurrences of the words in the
corpus of texts used to compile the dictionary.

Let us now consider books in a bookstore and let us consider the
entire collection of books sold. Assume that each book has an
inventory number. It is each copy of the book sold which is a
sign, and its value (price) is the cardinality $\omega$ of this
sign. The additional money involved in the value of the book
(storage expenses, overhead, etc.) must be included to get the
user cardinality $\widetilde{\omega}$ of the book. Here the
generator is the group of accountants who determined the prices
of books.

Now consider the catalog of books for sale. Each opus (be it a
novel, a collection of poems, or a textbook) is a sign, and its
price is the cardinality $\omega$ of the given sign.  Valuation
by user of this sign is the user  cardinality
$\widetilde{\omega}$. In this situation, all the sold copies of
the same opus, as opposed to the previous example, are grouped
together under one sign, which is specified by the title listed
in the catalog (this notion is similar to that of descriptor in
linguistics).

Each article of law in a book of statutes is a sign. The entire
list of laws is the dictionary of signs. Let us note that, in
specific examples, one can incorrectly interpret the notion of
sign and its cardinality. For instance, in the given example, the
number of people who were arrested under the given article of the
law is not, as one might think, the real cardinality of this sign
(article), and the number of people who actually broke this
article of law (whether they were arrested or not) is not its user
cardinality. The person (or persons) who created the book of
statutes is not the generator of the signs. Similarly, in
linguistics, the word forms that constitute a lexeme is not its
cardinality. Word forms are actually signs in a lower hierarchy.

The real cardinality  ${\omega}$ of articles of law, regarded as
signs, is the corresponding fine or the length of the prison
term. The user cardinality  $\widetilde{\omega}$ also includes all
the unpleasantness related to the punishment (the quality of
one's CV, separation from relatives, etc.). Here the generators
are the lawmakers who specified the punishment for breaking the
law.

People sent to prison for breaking the law may also be regarded
as signs (the inmates are even given serial numbers). The
cardinality  ${\omega}$ is the length of the prison term. The
generators in this case are the lawmakers, the judges, the
prosecuting attorneys, etc.

The sale of various goods will be considered below. Each type of
goods will be a sign, its price is its cardinality. The generator
is the person who fixed the prices. The set of all purchased
types of goods is the dictionary of signs, the prices are the
cardinalities, the customer is the user.

In the last two examples, it is easy to confuse the sign with its
cardinality. The user cardinality in these examples is quite
realistic for the customers, say for those who are buying cars.
Thus to obtain the user cardinality (price) of an automobile, one
must add to its list price the actual expenses related to its
upkeep, storage, insurance, spare parts, etc.

Similarly, in the case of judicial punishment, the cardinality
related to the actual losses for the prisoner (a spoiled CV, the
alienation by the family, etc.) becomes quite real for the other
people involved: for some the prisoner becomes an outcast, in
some cases becomes a hero for others.

In such cases, the user cardinality is not a monotone function of
the real cardinality. For cheap cars, it increases with the
decrease of the real cardinality, for instance, the insecurity of
the car becomes greater when its price decreases. Similarly, with
the decrease of prison terms, beginning at some level, the
related negative consequences do not decrease, and in fact
increase relatively to the real cardinality.

The generator should take into consideration the priorities, the
tastes, the possibilities, and so on of the user. If the
generator does not do this to a sufficient extent, then the
experimental curve will not approximate the theoretical curve as
well as it does in the automobile example shown  bellow on Figures
1 and 2.

For instance, if the generator (the lawmakers) does not take into
consideration the mentality of the given "user" and compiles a
set of laws under which practically any citizen constantly breaks
the laws, and, since it is impossible to imprison everyone, the
system starts putting in jail only those citizens which are in
power dislike for some reason, this will lead to a totalitarian
state where everyone lives in fear. In this case, the
experimental curves will not fit the theoretical ones, because
the absence of the preference principle (see \cite{NoPredp}) on
which the theoretical curves are based no longer applies.

Let us pass to the description of our main approach to the
general class of semiotic objets.

The most important and difficult question is how the generator
works out the cardinality of the dictionary of signs. These
cardinalities are worked out via a system of "agreements" between
the generator and the user. The generator "produces a fictional
action which places him at a higher level as compared to" the
user\footnote{See A. Grames, J. Courtier, Semiotics, An
Explanatory Dictionary, in \cite{Semiotica}.}.

If the generator, having recently passed the bar exam, begins to
impose an ideal system of laws to the user, it will be rejected
because it does not satisfy the social "rules of the game," and
the user will start reimplementing the lynching laws, or using
the laws of the maffia.

Let us look at another, even more spontaneous, generator, which
must include a huge number of people: it is impossible to specify
how and by whom the cities and towns of a country were founded.
What was the role of the interaction with neighbors, the greater
security in numbers, the role of commerce, all these factors must
be included in a very complicated and long algorithm.

Our considerations are based on Kolmogorov's approach to
randomness as maximal complexity (now known as Kolmogorov
complexity, see \cite{Kolmog1957}). This means that the longer
the algorithm used by the generator to construct the collection
of signs and their cardinalities, the nearer will the result be
to the general position of the majority of all possible versions
of these collections. This is similar to the fact that in playing
"heads or tails," the longer the number of trials, the nearer
will the sequence of heads and tails be to the "generic" version,
in which in half the trials we get heads, and tails in the other
half. And for the most part of the possible strings, we can apply
the theorem  from \cite{NoPredp}.

Indeed, let $N_i$ be the number of signs of the same real
cardinality ${\omega}$, while $\widetilde{\omega}$ is the user
cardinality. We denote the whole user
\textit{energia}\footnote{We use here  the  terminology of
Humbolt-Prieto}, by
$$
\cE= \sum_{i=1}^s N_i \widetilde{\omega_i}.
$$
We can assume that the number of signs $N_i$ corresponding to the
given user cardinality $\widetilde{\omega}$ of the sign $s_i$ is a
random variable with equiprobable distribution for any collection
of $\{N_i\}$ satisfying

\textbf{1)}
$$
\sum_{i=1}^s N_i \widetilde{\omega_i} \leq\cE,
$$
if
$$
\cE < \frac{\sum_{i=1}^s \widetilde{\omega_i}}{s} N;
$$
and \textbf{2)}
$$
\sum_{i=1}^s N_i \widetilde{\omega_i} \geq \cE,
$$
if
$$
\frac Ns \sum_{i=1}^s \widetilde{\omega_i}\leq
\cE\leq\widetilde{\omega}_{\max} N.
$$
Obviously, $\cE \leq \widetilde{\omega}_{\max} N$, where $N$ is
the length of the dictionary of signs.

This axiom should be understood in the sense that the given
string of signs of the energia $\cE$ is one of many such strings
with energia not greater than $\cE$, possessing the same
dictionary of signs; here we assume that, at least for the most
part of the signs, the energia is in general position with
respect to all possible versions of the collection $\{N_i\}$,
provided the latter satisfies conditions 1) or 2).

The case 1) has been proofed in \cite{DistrProb}. We present
bellow the proof of the case 2).

As in \cite{NegDimen}, the values of the random variable
$\widetilde{\omega}_1, \dots, \widetilde{\omega}_s$ are ordered
in absolute value.  In our consideration, both the number of
trials $N$ and~$s$ tend to infinity.

Let $N_i$ be the number of ''appearances'' of the value
$\widetilde{\omega}_i: \ \widetilde{\omega}_i \leq
\widetilde{\omega}_{i+1}$, then
\begin{equation}
\sum^s_{i=1} \frac{N_i}{N} \widetilde{\omega}_i=M, \label{Zipf1}
\end{equation}
where $M$ is the mathematical expectation.

The cumulative probability $\cP_k$ is the sum of the first~$k$
probabilities in the sequence $\widetilde{\omega}_i$:
$\cP_k=\frac 1N \sum_{i=1}^k N_i$, where $k<s$. We denote
$NP_k=B_k$.

If all the variants for which
\begin{equation}\label{A}
\sum_{i=1}^s N_i = N
\end{equation}
and
\begin{equation}\label{B}
\sum_{i=1}^s N_i \widetilde{\omega}_i \geq \cE, \ \ \cE=MN\ > N
\overline{\widetilde{\omega}},
\end{equation}
where $\overline{\widetilde{\omega}}=\frac{\sum_{i=1}^s
\widetilde{\omega}_i}{s}$, are equivalent (equiprobable), then
\cite{NoPredp} the majority of the variants will accumulate near
the following dependence of the ''cumulative probability''
$B_l\{N_i\}=\sum_{i=1}^l N_i$,
\begin{equation}
\sum_{i=1}^l N_i= \sum_{i=1}^l
\frac{1}{e^{\beta'\widetilde{\omega}_i-\nu'}-1}, \label{Zipf2}
\end{equation}
where $\beta'$ and $\nu'$ are determined by the conditions
\begin{equation}\label{Zipf2a}
B_s=N,
\end{equation}
\begin{equation}\label{Zipf2a'}
\sum_{i=1}^s \frac{ \widetilde{\omega}_i}{e^{\beta'
\widetilde{\omega}_i-\nu'}-1}=\cE,
\end{equation}
as $N \to \infty$ and $s \sim N$. By  the condition ~(\ref{B})
$\beta' <0$.

We introduce the notation: $\cM$ is the set of all sets $\{N_i\}$
satisfying conditions~(\ref{A}) and~(\ref{B}); \ $\cN\{\cM\}$ is
the number of elements of the set~$\cM$.

\begin{theorem} \label{theor1}
Suppose that all the variants of sets $\{N_i\}$ satisfying the
conditions ~(\ref{A}) and ~(\ref{B}) are equiprobable. Then the
number of variants $\cN$ of sets $\{N_i\}$ satisfying
conditions~(\ref{A}) and~(\ref{B}) and the additional relation
\begin{equation} \label{theorema1}
|\sum^l_{i=1} N_i - \sum^l_1\frac{1}{e^{\beta'
\widetilde{\omega}_i-\nu'}-1}|\geq  N^{(3/4+\varepsilon)}
\end{equation}
is less than $\frac{c_1 \cN\{\cM\}}{N^m}$ (where~$c_1$ and~$m$
are any arbitrary numbers, $l \geq\varepsilon N$, and
$\varepsilon$ is arbitrarily small).
\end{theorem}

{\bf{\qquad Proof of Theorem 1.}}

Let $\cA$ be a subset of $\cM$ satisfying the condition
$$
|\sum_{i=l+1}^s N_i - \sum_{i=l+1}^s \frac{1}
{e^{\beta\widetilde{\omega}_i-\nu}-1}|\leq \Delta;
$$
$$
|\sum_{i=1}^l N_i-\sum_{i=1}^l \frac{1}
{e^{\beta'\widetilde{\omega}_i-\nu'}-1}|\leq \Delta,
$$
where $\Delta$, $\beta$, $\nu$ are some real numbers independent
of~$l$.

We denote
$$
|\sum_{i=l+1}^s N_i-\sum_{i=l+1}^s \frac{1}
{e^{\beta\widetilde{\omega}_i-\nu}-1}| =S_{s-l};
$$
$$
|\sum_{i=1}^l N_i-\sum_{i=1}^l \frac{1}
{e^{\beta'\widetilde{\omega}_i-\nu'}-1}| =S_l.
$$

Obviously, if $\{N_i\}$ is the set of all sets of integers on the
whole, then
\begin{equation}\label{Proof1}
\cN\{\cM \setminus \cA\} = \sum_{\{N_i\}} \Bigl(
\Theta(\sum_{i=1}^s N_i\widetilde{\omega}_i -\cE)
\delta_{(\sum_{i=1}^s N_i),N} \Theta(S_l-\Delta)\Theta
(S_{s-l}-\Delta)\Bigr),
\end{equation}
where $\sum N_i=N$.

Here the sum is taken over all integers $N_i$,
$\Theta(\widetilde{\omega})$ is the Heaviside function, and
$\delta_{k_1,k_2}$ is the Kronecker symbol.

We use the integral representations
\begin{eqnarray}
&&\delta_{NN'}=\frac{e^{-\nu N}}{2\pi}\int_{-\pi}^\pi d\varphi
e^{-iN\varphi} e^{\nu N'}e^{i N'\varphi},\label{D7}\\
&&\Theta(y)=\frac1{2\pi i}\int_{-\infty}^\infty
d\widetilde{\omega}\frac1{\widetilde{\omega}-i}e^{\beta
y(1+i\widetilde{\omega})}.\label{D8}
\end{eqnarray}

Now we perform the standard regularization. We replace the first
Heaviside function~$\Theta$ in~(\ref{Proof1}) by the continuous
function
\begin{equation}
\Theta_{\alpha}(y) =\left\{
\begin{array}{ccc}
0 &\mbox{for}& \alpha > 1, \  y<0 \nn \\
1-e^{\beta y(1-\alpha)} &\mbox{for}& \alpha > 1,\  y \geq 0,
\label{Naz1}
\end{array}\right.
\end{equation}
\begin{equation}
\Theta_{\alpha}(y) =\left\{
\begin{array}{ccc}
e^{\beta y(1-\alpha)} &\mbox{for}&\alpha < 0, \ y<0 \nn \\
1 &\mbox{for}& \alpha < 0, \ y \geq 0, \label{Naz2}
\end{array}\right.
\end{equation}
where $\alpha \in (-\infty,0) \cup (1, \infty)$ is a parameter,
and obtain
\begin{equation}\label{proof2}
\Theta_\alpha(y) = \frac1{2\pi i} \int_{-\infty}^{\infty}
e^{\beta y(1+ix)} (\frac 1{x-i} - \frac 1{x-\alpha i}) dx.
\end{equation}

If  $\alpha > 1$, then $\Theta(y)\leq \Theta_{\alpha}(y)$.

Let $\nu <0$. We substitute~(\ref{D7}) and~(\ref{D8})
into~(\ref{Proof1}),  interchange the integration and summation,
then pass to the limit as $\alpha \to \infty$ and obtain the
estimate
\begin{eqnarray}
&&\cN\{\cM \setminus   \cA\} \leq \nn \\
&&\leq \Bigl|\frac{e^{-\nu N+\beta \cE
}}{i(2\pi)^2}\int_{-\pi}^\pi \bigl[ \exp(-iN\varphi)
\sum_{\{N_j\}}\exp\bigl\{-\beta\sum_{j=1}^s
N_j\widetilde{\omega}_j+(i\varphi+\nu)
\sum_{j=1}^s N_j\bigr\}\bigr]\ d\varphi \times \nn \\
&& \times \Theta(S_l -\Delta)\Theta(S_{s-l}-\Delta)\Bigr|,
\end{eqnarray}
where $\beta$ and $\nu$ are real parameters such that the series
converges for them.

To estimate the expression in the right-hand side, we bring the
absolute value sign inside the integral sign and then inside the
sum sign, integrate over $\varphi$, and obtain
\begin{eqnarray}
&&\cN\{\cM \setminus \cA\} \leq \frac{e^{-\nu N+\beta E }}{2\pi}
\sum_{\{N_i\}}\exp\{-\beta\sum_{i=1}^sN_i\widetilde{\omega}_i+\nu
\sum_{i=1}^s N_i\}\times \nn \\
&& \times\Theta (S_l-\Delta)\Theta (S_{s-l}-\Delta).
\end{eqnarray}

We denote
\begin{equation}\label{D9a}
Z(\beta,N)=\sum_{\{N_i\}} e^{-\beta\sum_{i=1}^s
N_i\widetilde{\omega}_i},
\end{equation}
where the sum is taken over all $N_i$ such that $\sum_{i=1}^s
N_i=N$,
$$
\zeta_l(\nu,\beta)= \prod_{i=1}^{l} \xi_i\left(\nu,\beta\right);
\zeta_{s-l}(\nu,\beta)= \prod_{i=l+1}^{s}
\xi_i\left(\nu,\beta\right);
$$
$$
\quad \xi_i(\nu,\beta)=
\frac{1}{(1-e^{\nu-\beta\widetilde{\omega}_i})}, \qquad
i=1,\dots,l.
$$

It follows from the inequality for the hyperbolic cosine
$\cosh(x)=(e^x+e^{-x})/2$ for $|x_1| \geq \delta; |x_2| \geq
\delta$:
\begin{equation}
\cosh(x_1)\cosh(x_2) > \frac{e^\delta}{2} \label{D33}
\end{equation}
that the inequality
\begin{equation}
\Theta(S_{s-l}-\Delta) \Theta(S_{l}-\Delta)\le e^{-c\Delta}
\cosh\Bigl(c\sum_{i=1}^{l} N_i
-c\phi_l\Bigr)\cosh\Bigl(c\sum_{i=l+1}^{s} N_i
-c\overline{\phi}_{s-l}\Bigr), \label{D34}
\end{equation}
where
$$
\phi_l= \sum_{i=1}^l
\frac{1}{e^{\beta'\widetilde{\omega}_i-\nu'}-1}; \qquad
\overline{\phi}_{s-l}= \sum_{i=l+1}^s
\frac{1}{e^{\beta\widetilde{\omega}_i-\nu}-1},
$$
holds for all positive $c$ and~$\Delta$.

We obtain
\begin{eqnarray}
&&\cN\{\cM \setminus \cA\} \leq  e^{-c\Delta} \exp\left(\beta
\cE-\nu N\right) \times \nn \\
&& \times \sum_{\{N_i\}}\exp\{-\beta\sum_{i=1}^l
N_i\widetilde{\omega}_i+\nu\sum_{i=1}^l N_i\}
\cosh\left(\sum_{i=1}^{l} c N_i -
c\phi\right) \times \nn \\
&& \times \exp\{-\beta\sum_{i=l+1}^s N_i\widetilde{\omega}_i +\nu
\sum_{i=l+1}^s
N_i\} \cosh\Bigl(\sum_{i=l+1}^s c N_i -c\overline{\phi}\Bigr) = \nn \\
&& =e^{\beta \cE} e^{-c\Delta} \times \nn \\
&&\times \left( \zeta_l(\nu-c,\beta) \exp(-c\phi_{l})
+\zeta_l(\nu+c,\beta)\exp(c\phi_{l})\right) \times \nn \\
&&\times\left(\zeta_{s-l}(\nu-c,\beta)
\exp(-c\overline{\phi}_{s-l})
+\zeta_{s-l}(\nu+c,\beta)\exp(c\overline{\phi}_{s-l})\right).
\label{5th}
\end{eqnarray}

Now we use the relations
\begin{equation}\label{5tha}
\frac {\pa}{\pa\nu}\ln \zeta_l|_{\beta=\beta',\nu=\nu'}\equiv
\phi_l; \frac {\pa}{\pa\nu}\ln
\zeta_{s-l}|_{\beta=\beta',\nu=\nu'}\equiv \overline{\phi}_{s-l}
\end{equation}
and the expansion $\zeta_l(\nu\pm c,\beta)$ by the Taylor formula.
There exists a $ \gamma <1$ such that
$$
\ln(\zeta_l(\nu\pm c,\beta)) =\ln\zeta_l(\nu,\beta)\pm
c(\ln\zeta_l)'_\nu(\nu,\beta)+\frac{c^2}{2}(\ln\zeta_l)^{''}_\nu
(\nu\pm\gamma c,\beta).
$$
We substitute this expansion, use formula~(\ref{5tha}), and see
that $\phi_{\nu,\beta}$ is cancelled.

Another representation of the Taylor formula implies
\begin{eqnarray}
&&\ln\left(\zeta_l(\nu+c,\beta)\right)=
\ln\left(\zeta_l(\beta,\nu)\right)+
\frac{c}\beta\frac{\pa}{\pa\nu}\ln\left(\zeta_l(\beta,\nu)\right)+\nn\\
&&+\int_{\nu}^{\nu+c/\beta}d\nu' (\nu+c/\beta-\nu')
\frac{\pa^2}{\pa\nu'^2}\ln\left(\zeta_l(\beta,\nu')\right).\label{CC1}
\end{eqnarray}
A similar expression holds for $\zeta_{s-l}$.

From the explicit form of the function $\zeta_l(\beta,\nu)$, we
obtain
\begin{equation}
\frac{\pa^2}{\pa\nu^2}\ln\left(\zeta_l(\beta,\nu)\right)=
\beta^2\sum_{i=1}^{l}
\frac{\exp(-\beta(\widetilde{\omega}_i+\nu))}{(\exp(-\beta(\widetilde{\omega}_i+\nu))-1)^2}\leq
\beta^2sd, \label{CC2}
\end{equation}
where $d$ is given by the formula
$$
d=\frac{\exp(-\beta(\widetilde{\omega}_s+\nu))}{(\exp(-\beta(\widetilde{\omega}_s+\nu))-1)^2}..
$$
The same estimate holds for $\zeta_{s-l}$.

Taking into account the fact that $\zeta_l\zeta_{s-l}=\zeta_s$,
we obtain the following estimate for $\beta=\beta'$ and
$\nu=\nu'$:
\begin{equation}\label{eval1}
\cN\{\cM \setminus \cA\}
\leq\zeta_s(\beta',\nu')\exp(-c\Delta+\frac{c^2}{2}\beta^2sd)
\exp(\cE\beta'-\nu'N).
\end{equation}

Now we express $\zeta_s(\nu',\beta')$ in terms $Z(\beta,N)$. To
do this, we prove the following lemma.

\begin{lemma}%lemma 1
Under the above assumptions, the asymptotics of the integral
\begin{equation}\label{lemma_1}
Z(\beta,N) = \frac{e^{-\nu N}}{2\pi}\int_{-\pi}^\pi
 d\alpha e^{-iN\alpha}\zeta_s(\beta,\nu+i\alpha)
\end{equation}
has the form
\begin{equation}\label{lemma_2}
Z(\beta,N) = C e^{-\nu N} \frac{\zeta_s(\beta,\nu)}{|(\partial^2
\ln\zeta_s(\beta,\nu))/ (\partial^2\nu)|} (1+O(\frac 1N)),
\end{equation}
where $C$ is a constant.
\end{lemma}

We have
\begin{equation}
 Z(\beta,N) = \frac{e^{-\nu N}}{2\pi}\int_{-\pi}^\pi
 e^{-iN\alpha}\zeta_s(\beta,\nu+i\alpha)\,d\alpha
 =\frac{e^{-\nu N}}{2\pi}\int_{-\pi}^\pi e^{NS(\alpha,N)} d\alpha ,\label{D15}
\end{equation}
where
\begin{equation}\label{qq}
    S(\alpha,N) = -i\alpha+ \ln \zeta_s (\beta, \nu +i\alpha)
    = -i\alpha - \sum_{i=1}^s \ln [1-e^{\nu+i\alpha-\beta\widetilde{\omega}_i}].
\end{equation}
Here $S$ depends on $N$, because $s$, $\widetilde{\omega}_i$, and
$\nu$ also depend on~$N$; the latter is chosen so that the point
$\alpha=0$ be a stationary point of the phase~$S$, i.e., from the
condition
\begin{equation}\label{qq1}
 N=\sum_{i=1}^s\frac{1}{e^{\beta\widetilde{\omega}_i-\nu}-1}.
\end{equation}
We assume that $a_1N \leq s \leq a_2N$, $a_1,
a_2=\operatorname{const}$, and, in addition,
$0\le\widetilde{\omega}_i\le B$ and $B=\operatorname{const}$,
$i=1,\dots,s$.
 If these conditions are satisfied
in some interval $\beta\in[0,\beta_0]$ of the values of the
inverse temperature, then all the derivatives of the phase are
bounded, the stationary point is nondegenerate, and the real part
of the phase outside a neighborhood of zero is strictly less than
its value at zero minus some positive number. Therefore,
calculating the asymptotics of the integral, we can replace the
interval of integration $[-\pi,\pi]$ by the interval $[-\e,\e]$.
In this integral, we perform the change of variable
\begin{equation}\label{qqq}
  z=\sqrt {S(0,N)-S(\alpha,N)}.
\end{equation}
This function is holomorphic in the disk $\abs{\alpha}\le\e$ in
the complex $\alpha$-plane and has a holomorphic inverse for a
sufficiently small~$\e$. As a result, we obtain
\begin{equation}\label{qqq1}
    \int_{-\e}^\e e^{NS(\alpha,N)}
    d\alpha=e^{NS(0,N)}\int_\gamma e^{-Nz^2}f(z)\,dz,
\end{equation}
where the path $\gamma$ in the complex $z$-plane is obtained from
the interval $[-\e,\e]$ by the change~\eqref{qqq} and
\begin{equation}\label{qqq3}
    f(z)=\left(\frac{\partial\sqrt {S(0,N)-S(\alpha,N)}}
    {\partial\alpha}\right)^{-1}\bigg|_{\alpha=\alpha(z)}.
\end{equation}
For a small~$\e$ the path $\gamma$ lies completely inside the
double sector $\re(z^2)>c(\re z)^2$ for some $c>0$; hence it can
be ``shifted'' to the real axis so that the integral does not
change up to terms that are exponentially small in~$N$. Thus,
with the above accuracy, we have
\begin{equation}\label{qqq4}
  Z(\beta,N) =  \frac{e^{-\nu N}}{2\pi}\int_{-\e}^\e e^{-Nz^2}f(z)\,dz.
\end{equation}
Since the variable $z$ is now real, we can assume that the
function $f(z)$ is finite (changing it outside the interval of
integration), extend the integral to the entire axis (which again
gives an exponentially small error), and then calculate the
asymptotic expansion of the integral expanding the integrand in
the Taylor series in~$z$ with a remainder. This justifies that
the saddle-point method can be applied to the above integral in
our case.

\begin{lemma}%lemma 2
The quantity
\begin{equation}\label{qqq5}
\frac{1}{\cN(\cM)} \sum_{\{N_i\}} e^{-\beta\sum_{i=1}^s
N_i\widetilde{\omega}_i},
\end{equation}
where $\sum N_i =N$ and $\widetilde{\omega}_iN_i>
\cE+N^{1/2+\varepsilon}$, tends to zero faster than $N^{-k}$ for
any $k$, $\varepsilon>0$.
\end{lemma}

We consider the point of minimum in $\beta$ of the right-hand
side of ~(\ref{5th}) with $\nu(\beta,N)$ satisfying the condition
$$
\sum \frac{1}{e^{\beta\widetilde{\omega}_i-\nu(\beta,N)}-1} =N.
$$
It is easy to see that it satisfies condition~(\ref{Zipf2a}). Now
we assume that the assumption of the lemma is not satisfied.

Then for $\sum N_i=N$,  $\sum \widetilde{\omega}_i N_i\geq
\cE+N^{1/2+\varepsilon}$, we have
$$
e^{\beta \cE}\sum_{\{N_i\}} e^{-\beta\sum_{i=1}^s
N_i\widetilde{\omega}_i}\geq e^{(N^{1/2}+\varepsilon)\beta}.
$$
Obviously, $\beta\ll \frac{1}{\sqrt{N}}$ provides a minimum
of~(\ref{5th}) if the assumptions of Lemma~1 are satisfied, which
contradicts the assumption that the minimum in~$\beta$ of the
right-hand side of~(\ref{5th}) is equal to~$\beta'$.

We set $c=\frac\Delta{N^{1+\alpha}}$ in formula~(\ref{eval1})
after the substitution~(\ref{lemma_2}); then it is easy to see
that the ratio
$$
\frac{\cN(\cM \setminus\cA)}{\cN(\cM)}\approx \frac 1{N^m},
$$
where $m$ is an arbitrary integer, holds for
$\Delta=N^{3/4+\varepsilon}$. The proof of the theorem is
complete.

We prove a cumulative formula in which the densities coincide in
shape with the Bose--Einstein distribution with negative
temperature. The difference consists also in that, instead of the
set $\widetilde{\omega}_n$ of random variables or eigenvalues of
the Hamiltonian operator, the formula contains some of their
averages over the cells. In view of our theorem, the
$\varepsilon_i$, which are averages of the energy
$\widetilde{\omega}_k$ at the $i$th cell, are nonlinear averages
in the sense of Kolmogorov~\cite{NelinSred}.

Let us number the signs constituting the dictionary in the order
of increase of their cardinality, beginning with the minimal
cardinality $\omega_{\min}$. The signs that have the same
cardinality are ordered arbitrarily. The number of each sign in
this ordering will be called its rank and denoted by $r$. If $l$
is the number of the signs of cardinality $\omega_l$ (beginning
from $\omega_{\min}$), then by $r_l$ we shall denote the number
of all signs with cardinality less than or equal to $\omega_l$.
By $r_{-l}$ we shall denote the number of all signs with
cardinality greater than $\omega_l$, so that $r_l+r_{-l}$ is the
total number of all signs.

Exactly as in the article \cite{RJMP_Lingvostat}, we see that the
rank $r_l$ of the signs, ordered by increasing cardinality,
satisfies relations (3), (5), and (8) appearing in
\cite{RJMP_Lingvostat}.

Let us set
$$
\widetilde{\omega}_i =\omega_i(1+\alpha\omega^\gamma_i
+\alpha^{-1}\omega_i^{-\gamma}),
$$
then as $\beta\ll 1$
$$
r_l=\frac{c_1}{1+\alpha\omega_l^\gamma}+c_2; \ \ \omega_l =(\frac
1\alpha \frac{r_l}{r_{-l}})^{1/\gamma}.
$$

Figures 1 and 2 show how well the generators of the prices of
American automobiles estimate the demand. The first plot shows
the dependence of the number of cars sold at a price equal to or
less than UJ on the price, the second one, the dependence of the
number of the car in the "dictionary of cars" (this number can be
regarded as the detailed make of the car) in the increasing order
of the car prices. The generators here are the people who
determined the price. The point of inflection of the graph
corresponds to the price level where the additional expenses are
minimal. We see that this point is practically the same on both
plots. This means that the "agreement" between the generator and
the user in this case  reaches the high level.

\begin{figure}[h] %'cars
%\begin{minipage}{0.5\linewidth}
%\centering\epsfig{figure= S0x0p0,width=\linewidth}
\centering\epsfig{figure=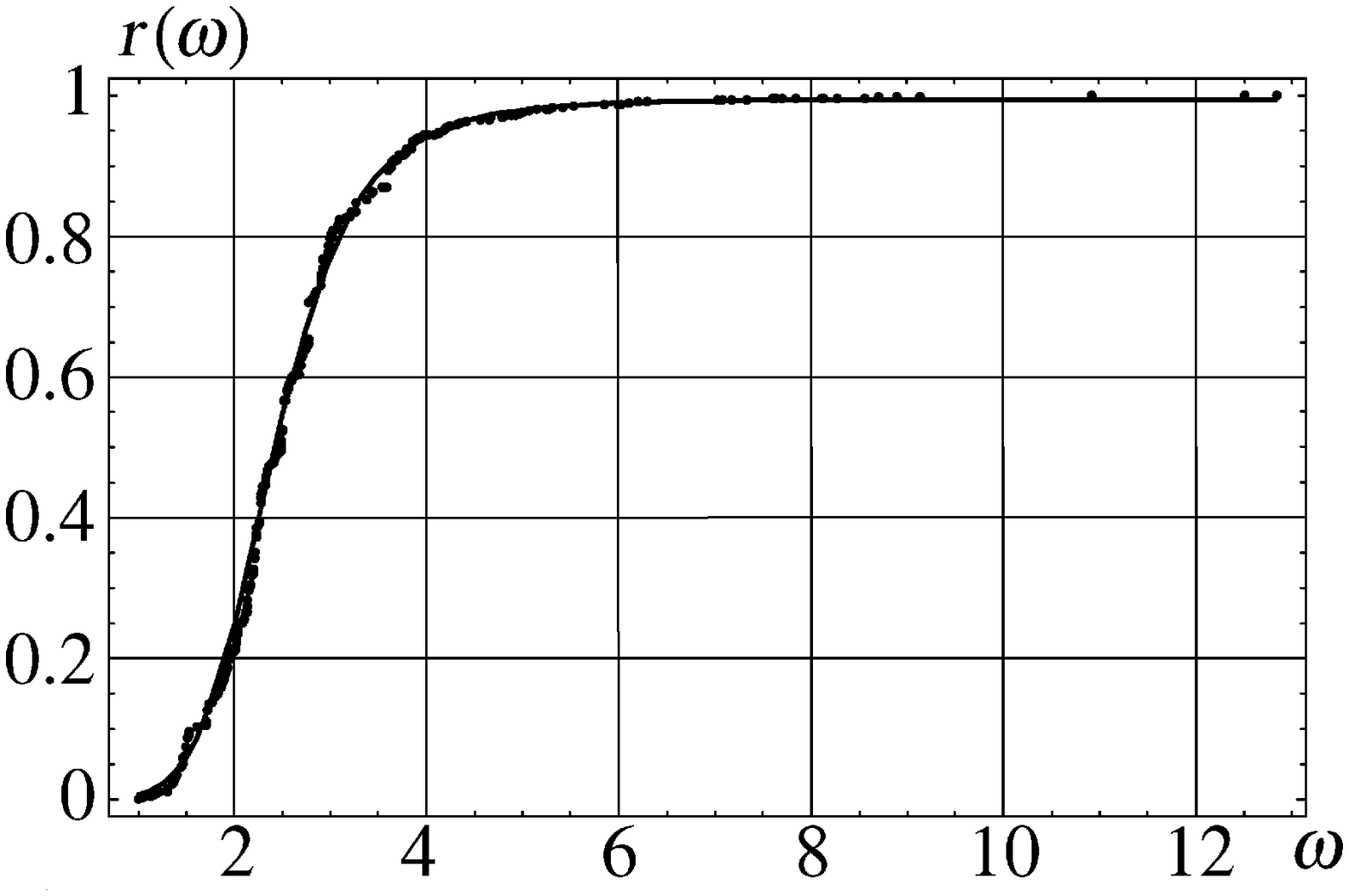,width=\linewidth} \caption{\small
Number of cars with price  $ < \ \omega$. The thin line represents
the theoretical curve $r(\omega)$. The mean quadratic error is
$\sigma = 0.0188674$.}.
%\end{minipage}\hfill
\end{figure}

\begin{figure}[h] %'cars
%\begin{minipage}{0.5\linewidth}
%\centering\epsfig{figure= S0x0p0,width=\linewidth}
\centering\epsfig{figure=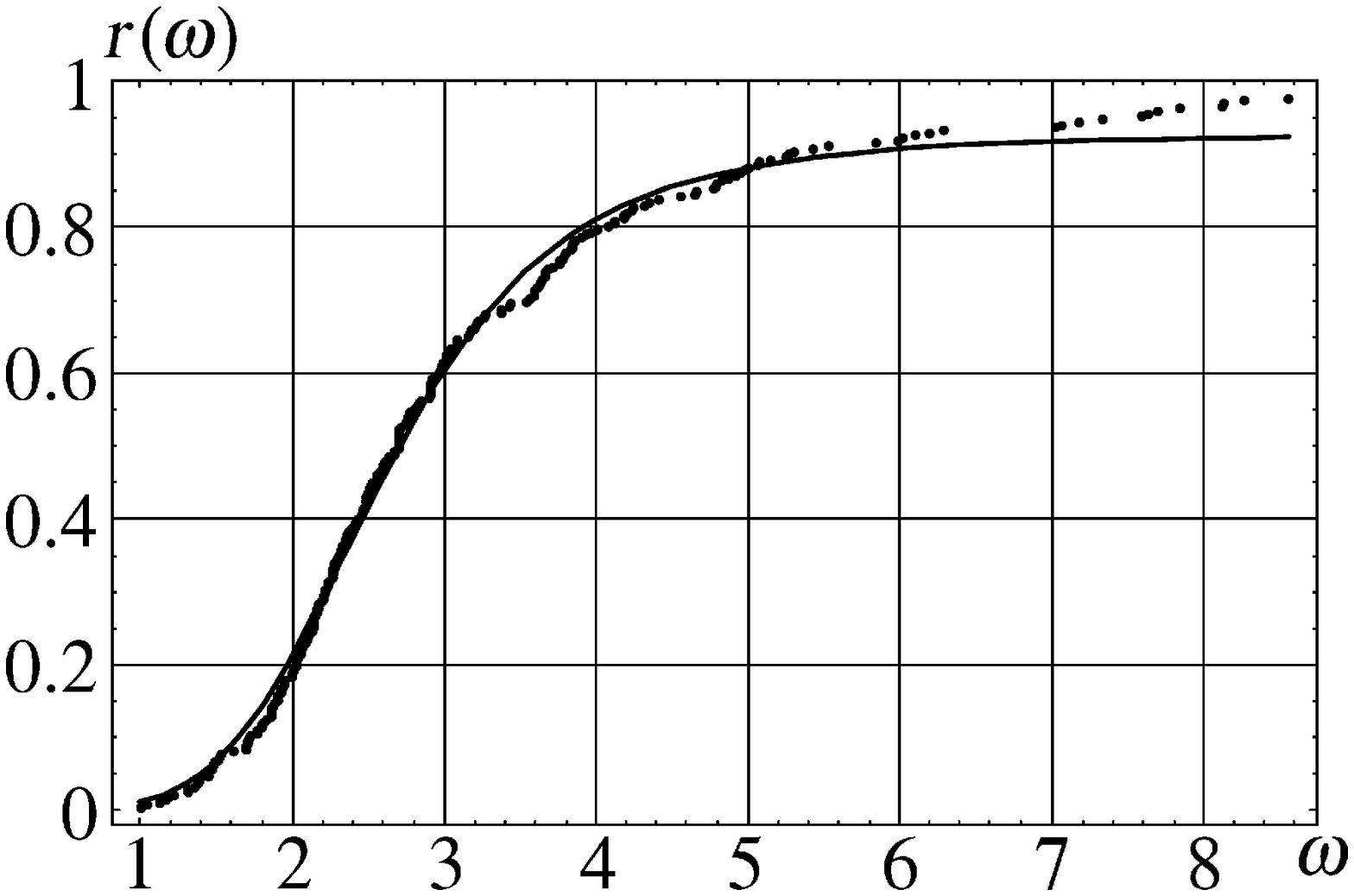,width=\linewidth} \caption{\small
Rank. The thin line represents the theoretical curve $r(\omega)$.
The mean quadratic error is $\sigma = 0.0184448$}.
%\end{minipage}\hfill
\end{figure}

\end{document}